\input amstex
\documentstyle{amsppt}
%----------------------------------------------------------------
% Title:     On Killing vector fields of a homogeneous and isotropic 
%            universe in closed model.
% Author:    Ruslan Sharipov
% Comments:  AmSTeX, 19 pages, amsppt style
% MSC-class: 53B30, 81T20, 83F05, 85A40
%----------------------------------------------------------------
%           Replacement for output macro definition
%
\catcode`@=11
\redefine\output@{%
  \def\break{\penalty-\@M}\let\par\endgraf
  \ifodd\pageno\global\hoffset=105pt\else\global\hoffset=8pt\fi  
  \shipout\vbox{%
    \ifplain@
      \let\makeheadline\relax \let\makefootline\relax
    \else
      \iffirstpage@ \global\firstpage@false
        \let\rightheadline\frheadline
        \let\leftheadline\flheadline
      \else
        \ifrunheads@ %\let\makefootline\relax
        \else \let\makeheadline\relax
        \fi
      \fi
    \fi
    \makeheadline \pagebody \makefootline}%
  \advancepageno \ifnum\outputpenalty>-\@MM\else\dosupereject\fi
}
\catcode`\@=\active
%----------------------------------------------------------------
\nopagenumbers
\def\negskp{\hskip -2pt}

\def\const{\operatorname{const}}

\def\blue#1{#1}
\catcode`#=11\def\diez{#}\catcode`#=6
\catcode`_=11\def\podcherkivanie{_}\catcode`_=8
%\catcode`~=11\def\volna{~}\catcode`~=\active
\def\mycite#1{\cite{\blue{#1}}\immediate\special{ps:
     ShrHPSdict begin /ShrBORDERthickness 0 def}}

\def\mytag#1{%
    \tag#1}
\def\mythetag#1{\thetag{\blue{#1}}\immediate\special{ps:
     ShrHPSdict begin /ShrBORDERthickness 0 def}}
\def\myrefno#1{\no#1}
\def\myhref#1#2{\blue{#2}\immediate\special{ps:
     ShrHPSdict begin /ShrBORDERthickness 0 def}}
\def\myEarXivlink{\myhref{http://arXiv.org}{http:/\negskp/arXiv.org}}
\def\myGeoCities{\myhref{http://www.geocities.com}{GeoCities}}
\def\mytheorem#1{\csname proclaim\endcsname{Theorem #1}}
\def\mythetheorem#1{\blue{#1}\immediate\special{ps:
     ShrHPSdict begin /ShrBORDERthickness 0 def}}
\def\mylemma#1{\csname proclaim\endcsname{Lemma #1}}

\def\mycorollary#1{\csname proclaim\endcsname{Corollary #1}}

%----------------------------------------------------------------
% Cyrillic fonts definition
%\font\eightcyr=wncyr8
%----------------------------------------------------------------
\pagewidth{360pt}
\pageheight{606pt}
\topmatter
\title
On Killing vector fields of a homogeneous and isotropic 
universe in closed model.
\endtitle
\author
R.~A.~Sharipov
\endauthor
\address 5 Rabochaya street, 450003 Ufa, Russia\newline
\vphantom{a}\kern 12pt Cell Phone: +7-(917)-476-93-48
\endaddress
\email \vtop to 30pt{\hsize=280pt\noindent
\myhref{mailto:r-sharipov\@mail.ru}
{r-sharipov\@mail.ru}\newline
\myhref{mailto:R\podcherkivanie Sharipov\@ic.bashedu.ru}
{R\_\hskip 1pt Sharipov\@ic.bashedu.ru}\vss}
\endemail
\urladdr
\vtop to 20pt{\hsize=280pt\noindent
\myhref{http://www.geocities.com/r-sharipov}
{http:/\negskp/www.geocities.com/r-sharipov}\newline
\myhref{http://www.freetextbooks.boom.ru/index.html}
{http:/\negskp/www.freetextbooks.boom.ru/index.html}\vss}
\endurladdr
\abstract
    Killing vector fields of a closed homogeneous and isotropic 
universe are studied. It is shown that in general case there is no
time-like Killing vector fields in such a universe. Two exceptional
cases are revealed.
\endabstract
\endtopmatter
\loadbold
\loadeufb
\TagsOnRight
\document
\accentedsymbol\tbvartheta{\tilde{\overline{\boldsymbol\vartheta}}
\vphantom{\boldsymbol\vartheta}}

\rightheadtext{On Killing vector fields \dots}
\head
1. Introduction. 
\endhead
\parshape 21 0pt 360pt 0pt 360pt 0pt 360pt 220pt 140pt 
220pt 140pt 220pt 140pt 220pt 140pt 220pt 140pt 220pt 140pt 
220pt 140pt 220pt 140pt 220pt 140pt 220pt 140pt 220pt 140pt 
220pt 140pt 220pt 140pt 220pt 140pt 220pt 140pt 220pt 140pt 
220pt 140pt 0pt 360pt
    Killing vector fields (infinitesimal isometries) are 
used in building vacuum states for quantum fields in a curved
space-time (see \mycite{1} and \mycite{2}). We study a homogeneous
and isotropic universe as an example of such a curved space-time.
\vadjust{\vskip 5pt\hbox to 0pt{\kern -10pt
\includegraphics{Iso01.eps}\hss}\vskip -5pt}This 
universe is diffeomorphic to the Cartesian product $M=\Bbb R
\times S^3$ (see \S\,111 and \S\,112 in \mycite{3}). Its spinor 
structure was studied in \mycite{4}. We use some technique from
\mycite{4} in present paper. In particular, we use the stereographic
coordinates $x^0,\,x^1,\,x^2,\,x^3$ and $y^0$, $y^1,\,y^2,\,y^3$
as two local charts co\-vering the whole universe. We call them 
{\it North Pole stereogra\-phic coordinates\/} and {\it South Pole 
stereographic coordinates\/} respectively. The domain of the North 
Pole stereographic coordinates\linebreak  $x^0,\,x^1,\,x^2,\,x^3$ 
is the whole sphe\-re $S^3$ except for one point, which
is called the North Pole. Similarly, the South Pole stereographic
coordinates are defined on the whole sphere $S^3$ except for the
diametrically opposite point, which is called the South Pole. Below
are the transition functions relating the North Pole and South Pole
stereographic coordinates in the intersection of their domains:
$$
\xalignat 2
	&\hskip -2em
\Vmatrix y^0\\y^1\\ y^2 \\ y^3\endVmatrix=\frac{1}{|x|^2}
\,\Vmatrix |x|^2\,x^0\\x^1\\ x^2 \\ x^3\endVmatrix,
&&\Vmatrix x^0\\x^1\\ x^2 \\ x^3\endVmatrix=\frac{1}{|y|^2}
\,\Vmatrix |y|^2\,y^0\\y^1\\ y^2 \\ y^3\endVmatrix.
\mytag{1.1}
\endxalignat
$$
Here $|x|^2=(x^1)^2+(x^2)^2+(x^3)^2$ and $|y|^2=(y^1)^2+(y^2)^2+(y^3)^2$.
The Minkowski type metric $\bold g$ in $M$ is given by the following
formulas:
$$
\align
&\hskip -2em
ds^{\kern 0.5pt 2}=R^2\,(dx^0)^2
-\frac{4\,R^2\,(dx^1)^2+4\,R^2\,(dx^2)^2
+4\,R^2\,(dx^3)^2}{\left(|x|^2+1\right)^2},
\mytag{1.2}\\
\vspace{2ex}
&\hskip -2em
ds^{\kern 0.5pt 2}=R^2\,(dy^0)^2
-\frac{4\,R^2\,(dy^1)^2+4\,R^2\,(dy^2)^2
+4\,R^2\,(dy^3)^2}{\left(|y|^2+1\right)^2}.
\mytag{1.3}
\endalign
$$
Looking at \mythetag{1.2} and \mythetag{1.3}, we see that the formulas
for metric in two different stereographic coordinates are very similar.
Therefore, we can derive some formulas in North Pole stereographic
coordinates and then transform them to South Pole coordinates by
substituting $y^0$, $y^1,\,y^2,\,y^3$ for $x^0,\,x^1,\,x^2,\,x^3$ 
without use of the transition functions \mythetag{1.1}.
    The parameter $R$ in the formulas \mythetag{1.2} and \mythetag{1.3}
is interpreted as the radius of the sphere $S^3$ in its realization as 
a hypersurface in the Euclidean space $\Bbb R^4$. This parameter is not 
a constant:
$$
\hskip -2em
R=R(x^0)=R(y^0).
\mytag{1.4}
$$
According to \mycite{3}, the time variable $t$ is introduced through 
the following formula:
$$
\hskip -2em
R\,dx^0=R\,dy^0=c\,dt\qquad\text{($c$ is the light velocity)}.
\mytag{1.5}
$$
Then we can write \mythetag{1.4} as $R=R(t)$. If $R(t)$ is constant,
we say that the universe is stable, if $R(t)$ is an increasing function,
we say that the universe is expanding, and if $R(t)$ is a decreasing 
function, we say that the universe is contracting. Oscillatory regimes 
are also possible. The main goal in this paper is to study under
which conditions for the function \mythetag{1.4} the universe
$M=\Bbb R\times S^3$ has at least one time-like Killing vector
field.
\head
2. Connection components and the curvature tensor.
\endhead
     The metric tensor $\bold g$ is determined by a diagonal matrix
$g_{ij}$ in North Pole stereographic coordinates. Its components are
determined by the formula \mythetag{1.2}:
$$
\xalignat 2
&\hskip -2em
g_{00}=R^2, &&g_{11}=g_{22}=g_{33}=-\frac{4\,R^2}{(|x|^2+1)^2}.
\mytag{2.1}
\endxalignat
$$
The dual metric tensor $\bold g$ is also given by a diagonal matrix:
$$
\xalignat 2
&\hskip -2em
g^{00}=\frac{1}{R^2}, 
&&g_{11}=g_{22}=g_{33}=-\frac{(|x|^2+1)^2}{4\,R^2}.
\mytag{2.2}
\endxalignat
$$
Now we choose the following well-known formula in order to 
calculate the components of the symmetric Levi-Civita connection:
$$
\hskip -2em
\Gamma^k_{ij}=\sum^3_{s=0}\frac{g^{ks}}{2}\left(\frac{\partial g_{is}}
{\partial x^j}+\frac{\partial g_{sj}}{\partial x^i}-\frac{\partial g_{ij}}
{\partial x^s}\right).
\mytag{2.3}
$$
Note that, unlike \mycite{4}, here we do not use non-holonomic frames
since we do not need to deal with spinors in this paper. Applying
\mythetag{2.3} to \mythetag{2.1} and \mythetag{2.2}, we get the
following complete list of the nonzero components of $\Gamma^k_{ij}$:
$$
\xalignat 3
&\hskip -2em
\Gamma^0_{11}=\frac{4\,R'}{R\,(|x|^2+1)^2},
&&\Gamma^0_{22}=\frac{4\,R'}{R\,(|x|^2+1)^2},
&&\Gamma^0_{33}=\frac{4\,R'}{R\,(|x|^2+1)^2},
\qquad\\
\vspace{1ex}
&\hskip -2em
\Gamma^1_{11}=-\frac{2\,(x^1)}{|x|^2+1},
&&\Gamma^1_{22}=\frac{2\,(x^1)}{|x|^2+1},
&&\Gamma^1_{33}=\frac{2\,(x^1)}{|x|^2+1},
\qquad\\
\vspace{-1ex}
&&&\mytag{2.4}\\
\vspace{-1ex}
&\hskip -2em
\Gamma^2_{11}=\frac{2\,(x^2)}{|x|^2+1},
&&\Gamma^2_{22}=-\frac{2\,(x^2)}{|x|^2+1},
&&\Gamma^2_{33}=\frac{2\,(x^2)}{|x|^2+1},
\qquad\\
\vspace{1ex}
&\hskip -2em
\Gamma^3_{11}=\frac{2\,(x^3)}{|x|^2+)},
&&\Gamma^3_{22}=\frac{2\,(x^3)}{|x|^2+1},
&&\Gamma^3_{33}=-\frac{2\,(x^3)}{|x|^2+1},
\qquad\\
\vspace{1ex}
\endxalignat
$$
%\vskip -4ex
$$
\xalignat 2
&\hskip -2em
\Gamma^2_{12}=\Gamma^2_{21}=-\frac{2\,(x^1)}{|x|^2+1},
&&\Gamma^3_{13}=\Gamma^3_{31}=-\frac{2\,(x^1)}{|x|^2+1},\\
\vspace{1ex}
&\hskip -2em
\Gamma^3_{23}=\Gamma^3_{32}=-\frac{2\,(x^2)}{|x|^2+1},
&&\Gamma^1_{21}=\Gamma^1_{12}=-\frac{2\,(x^2)}{|x|^2+1},\\
\vspace{1ex}
&\hskip -2em
\Gamma^1_{31}=\Gamma^1_{13}=-\frac{2\,(x^3)}{|x|^2+1},
&&\Gamma^2_{32}=\Gamma^2_{23}=-\frac{2\,(x^3)}{|x|^2+1},
\qquad
\mytag{2.5}\\
\vspace{1ex}
&\hskip -2em
\Gamma^0_{00}=\frac{R'}{R},
&&\Gamma^1_{01}=\Gamma^1_{10}=\frac{R'}{R},\\
\vspace{1ex}
&\hskip -2em
\Gamma^2_{02}=\Gamma^2_{20}=\frac{R'}{R},
&&\Gamma^3_{03}=\Gamma^3_{30}=\frac{R'}{R}.
\endxalignat
$$
Note that here we have 31 nonzero connection components, while in
\mycite{4} in the case of a non-holonomic frame we had 18.\par
      The next step is to calculate the components of the Riemannian
curvature tensor $\bold R$. They are given by the well-known formula
$$
\hskip -2em
R^p_{qij}=\frac{\partial\Gamma^p_{\!j\,q}}{\partial x^i}
-\frac{\partial\Gamma^p_{\!i\,q}}{\partial x^j}
+\sum^3_{h=0}\left(\Gamma^p_{\!i\,h}\,\Gamma^h_{\!j\,q}
-\Gamma^p_{\!j\,h}\,\Gamma^h_{\!i\,q}\right).
\mytag{2.6}
$$
Applying \mythetag{2.6} to \mythetag{2.4} and \mythetag{2.5}, we
derive the following expressions 
$$
\allowdisplaybreaks
\xalignat 2
&\hskip -1ex
R^0_{101}=-R^0_{110}=4\,\frac{R''\,R-(R')^2}{R^2\,(|x|^2+1)^2},
&&R^1_{001}=-R^1_{010}=\frac{R''\,R-(R')^2}{R^2},	
\qquad\\
\vspace{2ex}
&\hskip -1ex
R^0_{202}=-R^0_{220}=4\,\frac{R''\,R-(R')^2}{R^2\,(|x|^2+1)^2},
&&R^2_{002}=-R^2_{020}=\frac{R''\,R-(R')^2}{R^2},
\qquad
\mytag{2.7}\\
\vspace{2ex}
&\hskip -1ex
R^0_{303}=-R^0_{330}=4\,\frac{R''\,R-(R')^2}{R^2\,(|x|^2+1)^2},
&&R^3_{003}=-R^3_{030}=\frac{R''\,R-(R')^2}{R^2},
\qquad
\endxalignat
$$
$$
\allowdisplaybreaks
\align
&\hskip -1ex
R^1_{212}=-R^1_{221}=R^2_{121}=-R^2_{112}
=4\,\frac{(R')^2+R^2}{R^2\,(|x|^2+1)^2},
\qquad\\
\vspace{1ex}
&\hskip -1ex
R^2_{323}=-R^2_{332}=R^3_{232}=-R^3_{223}
=4\,\frac{(R')^2+R^2}{R^2\,(|x|^2+1)^2},
\qquad
\mytag{2.8}\\
\vspace{1ex}
&\hskip -1ex
R^3_{131}=-R^3_{113}=R^1_{313}=-R^1_{331}
=4\,\frac{(R')^2+R^2}{R^2\,(|x|^2+1)^2}.
\qquad
\vspace{1ex}
\endalign
$$
Using the above formulas \mythetag{2.7} and \mythetag{2.8}, we can 
calculate the components of the Ricci tensor. Here is the list of 
nonzero ones of them:
$$
\xalignat 2
&\hskip -2em
R_{00}=\frac{3\,(R')^2-3\,R\,R''}{R^2},
&&R_{11}=\frac{8\,R^2+4\,(R')^2+4\,R\,R''}{R^2\,(|x|^2+1)^2},
\qquad\\
\vspace{-1ex}
&&&\mytag{2.9}\\
\vspace{-1ex}
&\hskip -2em
R_{22}=\frac{8\,R^2+4\,(R')^2+4\,R\,R''}{R^2\,(|x|^2+1)^2},
&&R_{33}=\frac{8\,R^2+4\,(R')^2+4\,R\,R''}{R^2\,(|x|^2+1)^2}.
\qquad
\endxalignat
$$
And finally, using \mythetag{2.9}, we calculate the scalar curvature:
$$
R_{\,\text{scalar}}=-\frac{6}{R^2}-\frac{6\,R''}{R^3}.
\mytag{2.10}
$$
As we see, the scalar curvature given by the formula \mythetag{2.10}
coincides with the scalar curvature calculated in \mycite{1} and 
\mycite{4} for this particular model of the universe.\par
\head
3. Differential equations for Killing vector fields.
\endhead
     Killing vector fields are also known as infinitesimal isometries.
Local one-parame\-tric diffeomorphism groups generated by these vector
fields are composed by isometries --- they preserve the metric 
$\bold g$. Therefore, if $\bold X$ is a Killing vector field in $M$,
then the Lie derivative $L_{\bold X}$, when applied to $\bold g$, 
yields zero:
$$
\hskip -2em
L_{\bold X}(\bold g)=0.
\mytag{3.1}
$$
In the coordinate form the equation \mythetag{3.1} is written as 
follows:
$$
\hskip -2em
\sum^3_{s=0}X^s\,\frac{\partial g_{ij}}{\partial x^s}
+\sum^3_{s=0}g_{sj}\frac{\partial X^s}{\partial x^i}
+\sum^3_{s=0}g_{is}\frac{\partial X^s}{\partial x^j}=0.
\mytag{3.2}
$$
Let's replace the partial derivatives in \mythetag{3.2} with the
covariant derivatives:
$$
\gather
\hskip -2em
\frac{\partial g_{ij}}{\partial x^s}=\nabla_{\!s}g_{ij}
+\sum^3_{k=0}\Gamma^k_{s\kern 0.5pt i}\,g_{kj}
+\sum^3_{k=0}\Gamma^k_{sj}\,g_{ik},
\mytag{3.3}\\
\hskip -2em
\frac{\partial X^s}{\partial x^i}=\nabla_{\!i}X^s
-\sum^3_{k=0}\Gamma^s_{i\kern 0.5pt k}\,X^k,
\mytag{3.4}\\
\hskip -2em
\frac{\partial X^s}{\partial x^j}=\nabla_{\kern -2pt j}X^s
-\sum^3_{k=0}\Gamma^s_{jk}\,X^k.
\mytag{3.5}
\endgather
$$
Substituting \mythetag{3.3}, \mythetag{3.4}, and \mythetag{3.5}
into \mythetag{3.2} and taking into account the symmetry of
$g_{ij}$ and $\Gamma^k_{ij}$ with respect to the indices
$i$ and $j$, we get
$$
\hskip -2em
\nabla_{\!s}g_{ij}+\sum^3_{s=0}g_{sj}\,\nabla_{\!i}X^s
+\sum^3_{s=0}g_{is}\,\nabla_{\kern -2pt j}X^s=0.
\mytag{3.6}
$$
Now let's remember that the metric $\bold g$ is concordant with its
metric connection $\Gamma$, i\.\,e\. $\nabla\bold g=0$. As a result
the equation \mythetag{3.6} is reduced to
$$
\hskip -2em
\sum^3_{s=0}g_{sj}\,\nabla_{\!i}X^s
+\sum^3_{s=0}g_{is}\,\nabla_{\kern -2pt j}X^s=0.
\mytag{3.7}
$$
In a metric manifold each vector field $\bold X$ is associated with
some unique covector field. This covector field is usually denoted 
by the same symbol $\bold X$. The components of such two associated 
vectorial and covectorial fields are related to each other through 
the index lowering and index raising procedures:
$$
\xalignat 2
&\hskip -2em
X_i=\sum^3_{j=0}g_{ij}\,X^j,
&&X^i=\sum^3_{j=0}g^{ij}\,X_j.
\mytag{3.8}
\endxalignat
$$
Applying \mythetag{3.8} to \mythetag{3.7} and taking into account
that $\nabla\bold g=0$, we derive
$$
\hskip -2em
\nabla_{\!i}X_j+\nabla_{\kern -2pt j}X_i=0.
\mytag{3.9}
$$\par
     The equation \mythetag{3.9} is a basic equation for Killing 
vector fields we are going to study in this paper. It is written 
in the covectorial form. Let's denote
$$
\hskip -2em
\nabla_{\!i}X_j=Y_{ij}\text{\ \ for \ }i<j.
\mytag{3.10}
$$
In terms of \mythetag{3.10} the equation \mythetag{3.9} can be 
rewritten in the following form:
$$
\hskip -2em
\nabla_{\!i}X_j=\cases\quad Y_{ij} &\text{for \ }i<j,\\
\quad 0 &\text{for \ }i=j,\\-Y_{j\kern 1pt i} &\text{for \ }i>j.
\endcases
\mytag{3.11}
$$
The equations \mythetag{3.11} look like a Pfaff system of first 
order PDE's if we treat $Y_{ij}$ as new undetermined functions.
However, in this case we need to write the differential equations 
for these functions. For this purpose let's differentiate
\mythetag{3.9}:
$$
\hskip -2em
\nabla_{\!k}\nabla_{\!i}X_j+\nabla_{\!k}\nabla_{\kern -2pt j}X_i=0.
\mytag{3.12}
$$
Then we triplicate the equations \mythetag{3.12} by means of the 
cyclic transposition of indices: $i\to j\to k\to i$. As a result 
we get other two copies of the equation \mythetag{3.12}:
$$
\align
&\hskip -2em
\nabla_{\!i}\nabla_{\kern -2pt j}X_k+\nabla_{\!i}\nabla_{\!k}X_j=0,
\mytag{3.13}\\
\vspace{1ex}
&\hskip -2em
\nabla_{\kern -2pt j}\nabla_{\!k}X_i+\nabla_{\kern -2pt j}
\nabla_{\!i}X_k=0.
\mytag{3.14}
\endalign
$$
Now let's add \mythetag{3.13} and \mythetag{3.14}, then subtract 
\mythetag{3.12} from them. As a result we get
$$
\hskip -2em
\nabla_{\!i}\nabla_{\kern -2pt j}X_k+\nabla_{\kern -2pt j}
\nabla_{\!i}X_k=
(\nabla_{\!k}\nabla_{\!i}-\nabla_{\!i}\nabla_{\!k})X_j+
(\nabla_{\!k}\nabla_{\kern -2pt j}-\nabla_{\kern -2pt j}
\nabla_{\!k})X_i.
\mytag{3.15}
$$
In order to transform the equality \mythetag{3.15} we use the
following well-known identity:
$$
\hskip -2em
(\nabla_{\!i}\nabla_{\kern -2pt j}-\nabla_{\kern -2pt j}\nabla_{\!i})
X_k=-\sum^3_{s=0}R^s_{kij}\,X_s.
\mytag{3.16}
$$
Here $R^s_{kij}$ are the components of the Riemannian curvature 
tensor $\bold R$ (see \mythetag{2.6}, \mythetag{2.7}, and \mythetag{2.8}
above). Applying \mythetag{3.16} to \mythetag{3.15}, we derive
$$
\hskip -2em
2\,\nabla_{\!i}\nabla_{\kern -2pt j}X_k=-\sum^3_{s=0}R^s_{k\kern 0.5pt ij}
\,X_s-\sum^3_{s=0}R^s_{j\kern 0.5pt ki}\,X_s-\sum^3_{s=0}R^s_{i\kern 1pt
kj}\,X_s.
\mytag{3.17}
$$
Now let's recall the following identities:
$$
\xalignat 2
&\hskip -2em
R^s_{ijk}+R^s_{ikj}=0,
&&R^s_{ijk}+R^s_{k\kern 0.5pt ij}+R^s_{j\kern 0.5pt ki}=0.
\mytag{3.18}
\endxalignat
$$
These are the well-known identities for the components of the curvature 
tensor. Their proof can be found in \mycite{5}. Applying \mythetag{3.18}
to \mythetag{3.17}, we get
$$
\hskip -2em
\nabla_{\!i}\nabla_{\kern -2pt j}X_k
=\sum^3_{s=0}R^s_{ij\kern 0.5pt k}\,X_s.
\mytag{3.19}
$$
If we remember the notations \mythetag{3.10}, then \mythetag{3.19} can
be rewritten as
$$
\hskip -2em
\nabla_{\!i}Y_{j\kern 0.5pt k}=\sum^3_{s=0}R^s_{ij\kern 0.5pt k}\,X_s.
\mytag{3.20}
$$
Both \mythetag{3.11} and \mythetag{3.20} form a complete system of
Pfaff equations for ten functions $X_0,\,X_1,\,X_2,\,X_3,\,Y_{01},\, 
Y_{02},\,Y_{03},\,Y_{12},\,Y_{13},\,Y_{23}$. The following theorem
is an immediate consequence of this observation.
\mytheorem{3.1} A four-dimensional space-time manifold $M$ can have
at most ten linearly independent Killing vector fields.
\endproclaim
The actual number of isometries depends on the so-called compatibility
conditions for the Pfaff equations \mythetag{3.11} and \mythetag{3.20}.
In order to derive these compatibility conditions, let's calculate 
$\nabla_{\!i}\nabla_{\kern -2pt j}X_k-\nabla_{\kern -2pt j}\nabla_{\!i}
X_k$ and $\nabla_{\!i}\nabla_{\kern -2pt j}Y_{p\kern 1pt q}
-\nabla_{\kern -2pt j}\nabla_{\!i}Y_{p\kern 1pt q}$ on the base of 
the Pfaff equations \mythetag{3.11} and \mythetag{3.20}. The inequalities
in \mythetag{3.11} produce many special cases that should be studied 
separately. In order to avoid this inconvenience we extend the definition
of $Y_{ij}$. In \mythetag{3.10} they are defined for $i<j$. Let's set
$$
\hskip -2em
Y_{ij}=\cases \quad 0 &\text{for \ }i=j;\\
-Y_{j\kern 0.5pt i}&\text{for \ }i>j.
\endcases
\mytag{3.21}
$$
Due to the extension \mythetag{3.21} of \mythetag{3.10} we can write 
\mythetag{3.11} as
$$
\hskip -2em
\nabla_{\!i}X_j=Y_{ij}
\mytag{3.22}
$$
for all $i$ and $j$, but we should keep in mind that only $6$ of 
$16$ components of the skew-symmetric tensorial field $\bold Y$ 
are independent. Due to the skew symmetry $R^s_{ij\kern 0.3pt k}
=-R^s_{i\kern 0.2pt kj}$ the extension \mythetag{3.21} of 
\mythetag{3.10} and the extension \mythetag{3.22} of \mythetag{3.11}
are compatible with \mythetag{3.20}. Now for $\nabla_{\!i}
\nabla_{\kern -2pt j}X_k-\nabla_{\kern -2pt j}\nabla_{\!i}X_k$ 
we have
$$
\hskip -2em
\nabla_{\!i}\nabla_{\kern -2pt j}X_k-\nabla_{\kern -2pt j}
\nabla_{\!i}X_k=\nabla_{\!i}Y_{j\kern 0.5pt k}
-\nabla_{\kern -2pt j}Y_{i\kern 0.5pt k}.
\mytag{3.23}
$$
Applying the equation \mythetag{3.20} to the right hand side of
\mythetag{3.23} and applying the identity \mythetag{3.16} to its 
left hand side, we derive
$$
\hskip -2em
-\sum^3_{s=0}R^s_{kij}\,X_s=\sum^3_{s=0}R^s_{ij\kern 0.2pt k}\,X_s
-\sum^3_{s=0}R^s_{j\kern 0.3pt i\kern 0.5pt k}\,X_s.
\mytag{3.24}
$$
It is easy to see that the compatibility condition \mythetag{3.24}
is fulfilled identically due to the properties \mythetag{3.18} of the 
curvature tensor $\bold R$.\par
     Now we proceed to the compatibility conditions
derived from $\nabla_{\!i}\nabla_{\kern -2pt j}Y_{p\kern 1pt q}
-\nabla_{\kern -2pt j}\nabla_{\!i}Y_{p\kern 1pt q}$. In this case,
applying the equation \mythetag{3.20}, we get 
$$
\hskip -2em
\nabla_{\!i}\nabla_{\kern -2pt j}Y_{p\kern 1pt q}
-\nabla_{\kern -2pt j}\nabla_{\!i}Y_{p\kern 1pt q}
=\sum^3_{s=0}\nabla_{\!i}(R^s_{jp\kern 0.5pt q}\,X_s)
-\sum^3_{s=0}\nabla_{\kern -2pt j}(R^s_{ip\kern 0.5pt q}\,X_s).
\mytag{3.25}
$$
The left hand side of the equality \mythetag{3.25} is transformed by 
means of the identity 
$$
\hskip -2em
\nabla_{\!i}\nabla_{\kern -2pt j}Y_{p\kern 1pt q}
-\nabla_{\kern -2pt j}\nabla_{\!i}Y_{p\kern 1pt q}
=-\sum^3_{s=0}R^s_{p\kern 0.5pt ij}\,Y_{s\kern 0.1pt q}
-\sum^3_{s=0}R^s_{q\kern 0.5pt ij}\,Y_{p\kern 0.1pt s}.
\mytag{3.26}
$$
The identity \mythetag{3.26} is a tensorial generalization of
\mythetag{3.16}. Applying \mythetag{3.26} to \mythetag{3.25}
and taking into account \mythetag{3.22}, we transform
\mythetag{3.25} as follows:
$$
\hskip -2em
\gathered
-\sum^3_{s=0}R^s_{p\kern 0.5pt ij}\,Y_{s\kern 0.1pt q}
-\sum^3_{s=0}R^s_{q\kern 0.5pt ij}\,Y_{p\kern 0.1pt s}
=\sum^3_{s=0}\nabla_{\!i}R^s_{jp\kern 0.5pt q}\,X_s\,+\\
+\sum^3_{s=0}R^s_{jp\kern 0.5pt q}\,Y_{i\kern 0.2pt s}
-\sum^3_{s=0}\nabla_{\kern -2pt j}R^s_{ip\kern 0.5pt q}\,X_s
-\sum^3_{s=0}R^s_{ip\kern 0.5pt q}\,Y_{js}.
\endgathered
\mytag{3.27}
$$
The equality \mythetag{3.27} is a non-trivial compatibility 
condition for the system of Pfaff equations \mythetag{3.22} 
and \mythetag{3.20}. In the next section we shall study this 
equality for our particular case, where $M=\Bbb R\times S^3$.
\par
\head
4. Simplifying the compatibility conditions.
\endhead
    Note that the compatibility equations \mythetag{3.27} contain
the covariant derivatives of the curvature tensor. Therefore, we
begin our study of \mythetag{3.27} with calculating these covariant
derivatives. They are given by the formula:
$$
\nabla_{\!s}R^p_{q\kern 0.5pt ij}=\frac{R^p_{q\kern 0.5pt ij}}
{\partial x^s}+\sum^3_{h=0}\Gamma^p_{sh}\,R^h_{q\kern 0.5pt ij}
-\sum^3_{h=0}\Gamma^h_{s\kern 0.5pt q}\,R^p_{h\kern 0.5pt ij}
-\sum^3_{h=0}\Gamma^h_{s\kern 0.5pt i}\,R^p_{q\kern 0.1pt hj}
-\sum^3_{h=0}\Gamma^h_{s\kern 0.5pt j}\,R^p_{q\kern 0.1pt ih}.
$$
We substitute \mythetag{2.4}, \mythetag{2.5}, \mythetag{2.7},
and \mythetag{2.8} into this formula and get the following list 
of nonzero components $\nabla_{\!s}R^p_{q\kern 0.5pt ij}$ in 
North Pole stereographic coordinates:
$$
\allowdisplaybreaks
\align
&\hskip -2em
\aligned
\nabla_{\!0}R^0_{101}=-\nabla_{\!0}R^0_{110}
&=\nabla_{\!0}R^0_{202}=\\
\vspace{1ex}
=-\nabla_{\!0}R^0_{220}&=\nabla_{\!0}R^0_{303}
=-\nabla_{\!0}R^0_{330}=\\
\vspace{1ex}
&=\frac{16\,(R')^3-20\,R''\,R'\,R+4\,R'''\,R^2}
{R^3\,(|x|^2+1)^2},
\endaligned
\mytag{4.1}\\
\vspace{2ex}
&\hskip -2em
\aligned
\nabla_{\!0}R^1_{001}=-\nabla_{\!0}R^1_{010}
&=\nabla_{\!0}R^2_{002}=\\
\vspace{1ex}
=-\nabla_{\!0}R^2_{020}&=\nabla_{\!0}R^3_{003}
=-\nabla_{\!0}R^3_{030}=\\
\vspace{1ex}
&=\frac{4\,(R')^3-5\,R''\,R'\,R+\,R'''\,R^2}
{R^3},\\
\endaligned
\mytag{4.2}\\
\vspace{3ex}
&\hskip -2em
\aligned
\nabla_{\!0}R^1_{212}=-\nabla_{\!0}R^1_{221}
&=\nabla_{\!0}R^2_{121}=-\nabla_{\!0}R^2_{112}=\\
\vspace{1ex}
\nabla_{\!0}R^2_{323}=-\nabla_{\!0}R^2_{332}
&=\nabla_{\!0}R^3_{232}=-\nabla_{\!0}R^3_{223}=\\
\vspace{1ex}
\nabla_{\!0}R^3_{131}=-\nabla_{\!0}R^3_{113}
&=\nabla_{\!0}R^1_{313}=-\nabla_{\!0}R^1_{331}=\\
\vspace{1ex}
&=\frac{-16\,(R')^3+8\,R''\,R'\,R-8\,R'\,R^2}
{R^3\,(|x|^2+1)^2},
\endaligned
\mytag{4.3}\\
\vspace{3ex}
&\hskip -2em
\aligned
\nabla_{\!1}R^0_{212}=-\nabla_{\!1}R^0_{221}
&=\nabla_{\!2}R^0_{323}=-\nabla_{\!2}R^0_{332}=\\
\vspace{1ex}
\nabla_{\!3}R^0_{131}=-\nabla_{\!3}R^0_{113}
&=\nabla_{\!1}R^0_{313}=-\nabla_{\!1}R^0_{331}=\\
\vspace{1ex}
\nabla_{\!2}R^0_{121}=-\nabla_{\!2}R^0_{112}
&=\nabla_{\!3}R^0_{232}=-\nabla_{\!3}R^0_{223}=\\
\vspace{1ex}
&=\frac{32\,(R')^3-16\,R''\,R'\,R+16\,R'\,R^2}
{R^3\,(|x|^2+1)^4},
\endaligned
\mytag{4.4}\\
\vspace{3ex}
&\hskip -2em
\aligned
\nabla_{\!1}R^2_{012}=-\nabla_{\!1}R^2_{021}
&=\nabla_{\!2}R^3_{023}=-\nabla_{\!2}R^3_{032}=\\
\vspace{1ex}
\nabla_{\!3}R^1_{031}=-\nabla_{\!3}R^1_{013}
&=\nabla_{\!1}R^3_{013}=-\nabla_{\!1}R^3_{031}=\\
\vspace{1ex}
\nabla_{\!2}R^1_{021}=-\nabla_{\!2}R^1_{012}
&=\nabla_{\!3}R^2_{032}=-\nabla_{\!3}R^2_{023}=\\
\vspace{1ex}
&=\frac{8\,(R')^3-4\,R''\,R'\,R+4\,R'\,R^2}
{R^3\,(|x|^2+1)^2},
\endaligned
\mytag{4.5}\\
\vspace{3ex}
&\hskip -2em
\aligned
\nabla_{\!1}R^1_{220}=-\nabla_{\!1}R^1_{202}
&=\nabla_{\!2}R^2_{330}=-\nabla_{\!2}R^2_{303}=\\
\vspace{1ex}
\nabla_{\!3}R^3_{110}=-\nabla_{\!3}R^3_{101}
&=\nabla_{\!1}R^1_{330}=-\nabla_{\!1}R^1_{303}=\\
\vspace{1ex}
\nabla_{\!2}R^2_{110}=-\nabla_{\!2}R^2_{101}
&=\nabla_{\!3}R^3_{220}=-\nabla_{\!3}R^3_{202}=\\
\vspace{1ex}
&=\frac{8\,(R')^3-4\,R''\,R'\,R+4\,R'\,R^2}
{R^3\,(|x|^2+1)^2},
\endaligned
\mytag{4.6}\\
\vspace{3ex}
&\hskip -2em
\aligned
\nabla_{\!1}R^2_{102}=-\nabla_{\!1}R^2_{120}
&=\nabla_{\!2}R^3_{203}=-\nabla_{\!2}R^3_{230}=\\
\vspace{1ex}
\nabla_{\!3}R^1_{301}=-\nabla_{\!3}R^1_{310}
&=\nabla_{\!1}R^3_{103}=-\nabla_{\!1}R^3_{130}=\\
\vspace{1ex}
\nabla_{\!2}R^1_{201}=-\nabla_{\!2}R^1_{210}
&=\nabla_{\!3}R^2_{302}=-\nabla_{\!3}R^2_{320}=\\
\vspace{1ex}
&=\frac{8\,(R')^3-4\,R''\,R'\,R+4\,R'\,R^2}
{R^3\,(|x|^2+1)^2}.
\endaligned
\mytag{4.7}
\endalign
$$
Now we substitute \mythetag{4.1}, \mythetag{4.2}, \mythetag{4.3},
\mythetag{4.4},  \mythetag{4.5}, \mythetag{4.6}, and  \mythetag{4.7}
into \mythetag{3.27}. We also substitute \mythetag{2.7} and  
\mythetag{2.8} into \mythetag{3.27}. As a result we obtain a series
of linear algebraic equations for the functions $X_0,\,X_1,\,X_2,
\,X_3,\,Y_{01},\,Y_{02},\,Y_{03},\,Y_{12},\,Y_{13},\,Y_{23}$.
Since the expressions in both sides of \mythetag{3.27} are 
skew-symmetric with respect to $i$ and $j$ and with respect to $p$ 
and $q$, we could have at most $36$ mutually independent equations.
However, in our particular case the number of mutually independent
equations is $5$. Here is the list of these five equations:
$$
\align
&\hskip -2em
\cases
\left(4\,(R')^3-5\,R''\,R'\,R+\,R'''\,R^2\right)X_0=0,\\
\vspace{1ex}
R'\left(2\,(R')^2-R''\,R+R^2\right)X_0=0,
\endcases
\mytag{4.8}\\
\vspace{2ex}
&\hskip -2em
\cases
\left(2\,(R')^2-R''\,R+R^2\right)\left(R'\,X_1-R\,Y_{01}\right)=0,\\
\left(2\,(R')^2-R''\,R+R^2\right)\left(R'\,X_2-R\,Y_{02}\right)=0,\\
\left(2\,(R')^2-R''\,R+R^2\right)\left(R'\,X_3-R\,Y_{03}\right)=0.
\endcases
\mytag{4.9}
\endalign
$$
As we see the compatibility equations \mythetag{4.8} and \mythetag{4.9}
depend essentially on the function \mythetag{1.4} and its derivatives.
These simplified equations will be studied in the next two sections.
\head
5. Spacial rotations.
\endhead
    Note that the functions $Y_{12},\,Y_{13},\,Y_{23}$ are not presented
in the equations \mythetag{4.8} and \mythetag{4.9}. This fact reflects
the spherical symmetry of our universe $M=\Bbb R\times S^3$. It is 
known that the sphere $S^3$ has a $6$-parametric group of isometries. 
These isometries produce $6$ linearly independent Killing vector fields
corresponding to $3$ meridional and $3$ equatorial rotations. Now we
write these vector fields explicitly.\par
    {\bf Meridional rotation} in the plane $z^1 O z^4$. This rotation 
induces a Killing vector field $\bold X$ in $M$ expressed by the formula
$$
\hskip -2em
\bold X=\frac{2\,(x^1)^2-|x|^2+1}{2}\,\frac{\partial}{\partial x^1}
+(x^1)\,(x^2)\,\frac{\partial}{\partial x^2}
+(x^1)\,(x^3)\,\frac{\partial}{\partial x^3}
\mytag{5.1}
$$
in the North Pole stereographic coordinates. Applying the index lowering
procedure \mythetag{3.8} to the components of \mythetag{5.1}, we get
the covectorial components of $\bold X$:
$$
\xalignat 2
&\hskip -2em
X_0=0, 
&&X_1=-\frac{4\,(x^1)^2-2\,|x|^2+2}{(|x|^2+1)^2}\,R^2,
\qquad\\
\vspace{-1ex}
\mytag{5.2}\\
\vspace{-1ex}
&\hskip -2em
X_2=-\frac{4\,(x^1)\,(x^2)}{(|x|^2+1)^2}\,R^2,
&&X_3=-\frac{4\,(x^1)\,(x^3)}{(|x|^2+1)^2}\,R^2.
\qquad
\endxalignat 
$$
Then, using the formula \mythetag{3.10}, we calculate the functions
$Y_{01},\,Y_{02},\,Y_{03},\,Y_{12},\,Y_{13}$, $Y_{23}$ associated 
with the Killing vector field \mythetag{5.1}: 
$$
\xalignat 2
&\hskip -2em
Y_{01}=-\frac{2\,R\,R'\,(2\,(x^1)^2-|x|^2+1)}{(|x|^2+1)^2},
&&Y_{02}=-\frac{4\,R\,R'\,(x^1)\,(x^2)}{(|x|^2+1)^2},\\
\vspace{1ex}
&\hskip -2em
Y_{03}=-\frac{4\,R\,R'\,(x^1)\,(x^3)}{(|x|^2+1)^2},
&&Y_{12}=-\frac{8\,R^2\,(x^2)}{(|x|^2+1)^3},
\qquad
\mytag{5.3}\\
\vspace{1ex}
&\hskip -2em
Y_{23}=0,
&&Y_{13}=-\frac{8\,R^2\,(x^3)}{(|x|^2+1)^3}.
\endxalignat 
$$
Let's substitute $x^1=x^2=x^3=0$ into \mythetag{5.2} and \mythetag{5.3}.
As a result we get
$$
\xalignat 3
&\hskip -2em
\Vmatrix X_0\\X_1\\X_2\\ X_3\endVmatrix
=\Vmatrix 0\\-2\,R^2\\ 0\\ 0\endVmatrix,
&&\Vmatrix Y_{01}\\Y_{02}\\Y_{03}\endVmatrix
=\Vmatrix -2\,R\,R'\\0\\ 0\endVmatrix,
&&\Vmatrix Y_{12}\\Y_{13}\\Y_{23}\endVmatrix
=\Vmatrix 0\\0\\ 0\endVmatrix.
\quad
\mytag{5.4}
\endxalignat
$$
The quantities listed in the formulas \mythetag{5.4} can be treated 
as initial data for the Pfaff equations \mythetag{3.20} and 
\mythetag{3.22}.\par
    {\bf Meridional rotation} in the plane $z^2 O z^4$. This case 
is very similar to the previous one. Here is the formula for the
Killing vector field in this case:
$$
\hskip -2em
\bold X=(x^2)\,(x^1)\,\frac{\partial}{\partial x^1}
+\frac{2\,(x^2)^2-|x|^2+1}{2}\,\frac{\partial}{\partial x^2}
+(x^2)\,(x^3)\,\frac{\partial}{\partial x^3}.
\mytag{5.5}
$$
Below are the covariant components of the vector \mythetag{5.5}
$$
\xalignat 2
&\hskip -2em
X_0=0, 
&&X_1=-\frac{4\,(x^2)\,(x^1)}{(|x|^2+1)^2}\,R^2,
\qquad\\
\vspace{-1ex}
\mytag{5.6}\\
\vspace{-1ex}
&\hskip -2em
X_2=-\frac{4\,(x^2)^2-2\,|x|^2+2}{(|x|^2+1)^2}\,R^2,
&&X_3=-\frac{4\,(x^2)\,(x^3)}{(|x|^2+1)^2}\,R^2.
\qquad
\endxalignat 
$$
Substituting \mythetag{5.6} into \mythetag{3.10} we obtain the
functions $Y_{01},\,Y_{02},\,Y_{03},\,Y_{12},\,Y_{13},\,Y_{23}$:
$$
\xalignat 2
&\hskip -2em
Y_{01}=-\frac{4\,R\,R'\,(x^2)\,(x^1)}{(|x|^2+1)^2},
&&Y_{02}=-\frac{2\,R\,R'\,(2\,(x^2)^2-|x|^2+1)}{(|x|^2+1)^2},\\
\vspace{1ex}
&\hskip -2em
Y_{03}=-\frac{4\,R\,R'\,(x^2)\,(x^3)}{(|x|^2+1)^2},
&&Y_{12}=\frac{8\,R^2\,(x^1)}{(|x|^2+1)^3},
\qquad
\mytag{5.7}\\
\vspace{1ex}
&\hskip -2em
Y_{23}=-\frac{8\,R^2\,(x^3)}{(|x|^2+1)^3},
&&Y_{13}=0.
\endxalignat 
$$
By setting $x^1=x^2=x^3=0$ in \mythetag{5.6} and \mythetag{5.7} we
derive
$$
\xalignat 3
&\hskip -2em
\Vmatrix X_0\\X_1\\X_2\\ X_3\endVmatrix
=\Vmatrix 0\\0\\-2\,R^2\\ 0\endVmatrix,
&&\Vmatrix Y_{01}\\Y_{02}\\Y_{03}\endVmatrix
=\Vmatrix 0\\-2\,R\,R'\\ 0\endVmatrix,
&&\Vmatrix Y_{12}\\Y_{13}\\Y_{23}\endVmatrix
=\Vmatrix 0\\0\\ 0\endVmatrix.
\quad
\mytag{5.8}
\endxalignat
$$
Like \mythetag{5.4}, the quantities \mythetag{5.8} are initial 
data for the equations \mythetag{3.20} and \mythetag{3.22}.
\par
     {\bf Meridional rotation} in the plane $z^3 O z^4$. This case
is also very similar to the previous cases. Here is the formula for 
the Killing vector field in this case:
$$
\hskip -2em
\bold X=(x^3)\,(x^1)\,\frac{\partial}{\partial x^1}
+(x^3)\,(x^2)\,\frac{\partial}{\partial x^2}
+\frac{2\,(x^3)^2-|x|^2+1}{2}\,\frac{\partial}{\partial x^3}.
\mytag{5.9}
$$
Below are the covariant components of the vector \mythetag{5.9}:
$$
\xalignat 2
&\hskip -2em
X_0=0, 
&&X_1=-\frac{4\,(x^3)\,(x^1)}{(|x|^2+1)^2}\,R^2,
\qquad\\
\vspace{-1ex}
\mytag{5.10}\\
\vspace{-1ex}
&\hskip -2em
X_2=-\frac{4\,(x^3)\,(x^2)}{(|x|^2+1)^2}\,R^2,
&&X_3=-\frac{4\,(x^3)^2-2\,|x|^2+2}{(|x|^2+1)^2}\,R^2.
\qquad
\endxalignat 
$$
Now, substituting \mythetag{5.10} into \mythetag{3.10}, we find
$Y_{01},\,Y_{02},\,Y_{03},\,Y_{12},\,Y_{13},\,Y_{23}$:
$$
\xalignat 2
&\hskip -2em
Y_{01}=-\frac{4\,R\,R'\,(x^3)\,(x^1)}{(|x|^2+1)^2},
&&Y_{02}=-\frac{4\,R\,R'\,(x^3)\,(x^2)}{(|x|^2+1)^2},\\
\vspace{1ex}
&\hskip -2em
Y_{03}=-\frac{2\,R\,R'\,(2\,(x^3)^2-|x|^2+1)}{(|x|^2+1)^2},
&&Y_{12}=0,
\qquad
\mytag{5.11}\\
\vspace{1ex}
&\hskip -2em
Y_{23}=\frac{8\,R^2\,(x^2)}{(|x|^2+1)^3},
&&Y_{13}=\frac{8\,R^2\,(x^1)}{(|x|^2+1)^3}.
\endxalignat 
$$
By setting $x^1=x^2=x^3=0$ in \mythetag{5.10} and \mythetag{5.11} we
obtain
$$
\xalignat 3
&\hskip -2em
\Vmatrix X_0\\X_1\\X_2\\ X_3\endVmatrix
=\Vmatrix 0\\0\\0\\-2\,R^2\endVmatrix,
&&\Vmatrix Y_{01}\\Y_{02}\\Y_{03}\endVmatrix
=\Vmatrix 0\\0\\-2\,R\,R'\endVmatrix,
&&\Vmatrix Y_{12}\\Y_{13}\\Y_{23}\endVmatrix
=\Vmatrix 0\\0\\ 0\endVmatrix.
\quad
\mytag{5.12}
\endxalignat
$$
The quantities \mythetag{5.12} are initial data for the equations
\mythetag{3.20} and \mythetag{3.22}.\par
    The next three cases are produced by the equatorial rotations.
They are somewhat different from meridional ones.\par
    {\bf Equatorial rotation} in the plane $z^1 O z^2$. The Killing 
vector field in this case is given by the following formula:
$$
\hskip -2em
\bold X=(x^2)\,\frac{\partial}{\partial x^1}
-(x^1)\,\frac{\partial}{\partial x^2}.
\mytag{5.13}
$$
This formula is more simple than \mythetag{5.1}, \mythetag{5.5}, or
\mythetag{5.9}. Here are the covariant components of the vector
field given by the formula \mythetag{5.13}:
$$
\xalignat 2
&\hskip -2em
X_0=0, 
&&X_1=-\frac{4\,R^2\,(x^2)}{(|x|^2+1)^2}\,R^2,
\qquad\\
\vspace{-1ex}
\mytag{5.14}\\
\vspace{-1ex}
&\hskip -2em
X_2=\frac{4\,R^2\,(x^1)}{(|x|^2+1)^2}\,R^2,
&&X_3=0.
\qquad
\endxalignat 
$$
Substituting \mythetag{5.14} into \mythetag{3.10}, we calculate the
functions $Y_{01},\,Y_{02},\,Y_{03},\,Y_{12},\,Y_{13}$, $Y_{23}$
for the vector field given by the formula \mythetag{5.13}:
$$
\xalignat 2
&\hskip -2em
Y_{01}=-\frac{4\,R\,R'\,(x^2)}{(|x|^2+1)^2},
&&Y_{02}=\frac{4\,R\,R'\,(x^1)}{(|x|^2+1)^2},\\
\vspace{1ex}
&\hskip -2em
Y_{03}=0,
&&Y_{12}=\frac{4\,R^2\,(2\,(x^3)^2-|x|^2+1)}{(|x|^2+1)^3},
\qquad
\mytag{5.15}\\
\vspace{1ex}
&\hskip -2em
Y_{23}=\frac{8\,R^2\,(x^3)\,(x^1)}{(|x|^2+1)^3},
&&Y_{13}=-\frac{8\,R^2\,(x^3)\,(x^2)}{(|x|^2+1)^3}.
\endxalignat 
$$
By setting $x^1=x^2=x^3=0$ in \mythetag{5.14} and \mythetag{5.15},
we obtain
$$
\xalignat 3
&\hskip -2em
\Vmatrix X_0\\X_1\\X_2\\ X_3\endVmatrix
=\Vmatrix 0\\0\\0\\0\endVmatrix,
&&\Vmatrix Y_{01}\\Y_{02}\\Y_{03}\endVmatrix
=\Vmatrix 0\\0\\0\endVmatrix,
&&\Vmatrix Y_{12}\\Y_{13}\\Y_{23}\endVmatrix
=\Vmatrix 4\,R^2\\0\\ 0\endVmatrix.
\quad
\mytag{5.16}
\endxalignat
$$
The quantities \mythetag{5.16} are initial data for the equations
\mythetag{3.20} and \mythetag{3.22}.\par
    {\bf Equatorial rotation} in the plane $z^2 O z^3$. This case 
is very similar to the previous one. Here is the formula for the
Killing vector field in this case:
$$
\hskip -2em
\bold X=(x^3)\,\frac{\partial}{\partial x^2}
-(x^2)\,\frac{\partial}{\partial x^3}.
\mytag{5.17}
$$
Below are the covariant components of the vector \mythetag{5.17}:
$$
\xalignat 2
&\hskip -2em
X_0=0, 
&&X_1=0,
\qquad\\
\vspace{-1ex}
\mytag{5.18}\\
\vspace{-1ex}
&\hskip -2em
X_2=-\frac{4\,R^2\,(x^3)}{(|x|^2+1)^2}\,R^2,
&&X_3=\frac{4\,R^2\,(x^2)}{(|x|^2+1)^2}\,R^2.
\qquad
\endxalignat 
$$
Substituting \mythetag{5.18} into \mythetag{3.10}, we calculate
$Y_{01},\,Y_{02},\,Y_{03},\,Y_{12},\,Y_{13},\,Y_{23}$:
$$
\xalignat 2
&\hskip -2em
Y_{01}=0,
&&Y_{02}=-\frac{4\,R\,R'\,(x^3)}{(|x|^2+1)^2},\\
\vspace{1ex}
&\hskip -2em
Y_{03}=\frac{4\,R\,R'\,(x^2)}{(|x|^2+1)^2},
&&Y_{12}=\frac{8\,R^2\,(x^1)\,(x^3)}{(|x|^2+1)^3},
\qquad
\mytag{5.19}\\
\vspace{1ex}
&\hskip -2em
Y_{23}=\frac{4\,R^2\,(2\,(x^1)^2-|x|^2+1)}{(|x|^2+1)^3},
&&Y_{13}=-\frac{8\,R^2\,(x^1)\,(x^2)}{(|x|^2+1)^3}.
\endxalignat 
$$
By setting $x^1=x^2=x^3=0$ in \mythetag{5.18} and \mythetag{5.19},
we obtain
$$
\xalignat 3
&\hskip -2em
\Vmatrix X_0\\X_1\\X_2\\ X_3\endVmatrix
=\Vmatrix 0\\0\\0\\0\endVmatrix,
&&\Vmatrix Y_{01}\\Y_{02}\\Y_{03}\endVmatrix
=\Vmatrix 0\\0\\0\endVmatrix,
&&\Vmatrix Y_{12}\\Y_{13}\\Y_{23}\endVmatrix
=\Vmatrix 0\\4\,R^2\\ 0\endVmatrix.
\quad
\mytag{5.20}
\endxalignat
$$
The quantities \mythetag{5.20} are initial data for the equations
\mythetag{3.20} and \mythetag{3.22}.\par
    {\bf Equatorial rotation} in the plane $z^3 O z^1$. The Killing 
vector field in this last case is given by the following formula:
$$
\hskip -2em
\bold X=-(x^3)\,\frac{\partial}{\partial x^1}
+(x^1)\,\frac{\partial}{\partial x^3}.
\mytag{5.21}
$$
Below are the covariant components of the vector \mythetag{5.21}:
$$
\xalignat 2
&\hskip -2em
X_0=0, 
&&X_1=\frac{4\,R^2\,(x^3)}{(|x|^2+1)^2}\,R^2,
\qquad\\
\vspace{-1ex}
\mytag{5.22}\\
\vspace{-1ex}
&\hskip -2em
X_2=0,
&&X_3=-\frac{4\,R^2\,(x^2)}{(|x|^2+1)^2}\,R^2.
\qquad
\endxalignat 
$$
Substituting \mythetag{5.22} into \mythetag{3.10}, we find
the functions $Y_{01},\,Y_{02},\,Y_{03},\,Y_{12},\,Y_{13},
\,Y_{23}$:
$$
\xalignat 2
&\hskip -2em
Y_{01}=\frac{4\,R\,R'\,(x^3)}{(|x|^2+1)^2},
&&Y_{02}=0,\\
\vspace{1ex}
&\hskip -2em
Y_{03}=-\frac{4\,R\,R'\,(x^1)}{(|x|^2+1)^2},
&&Y_{12}=\frac{8\,R^2\,(x^2)\,(x^3)}{(|x|^2+1)^3},
\qquad
\mytag{5.23}\\
\vspace{1ex}
&\hskip -2em
Y_{23}=-\frac{8\,R^2\,(x^2)\,(x^1)}{(|x|^2+1)^3},
&&Y_{13}=-\frac{4\,R^2\,(2\,(x^2)^2-|x|^2+1)}{(|x|^2+1)^3}.
\endxalignat 
$$
By setting $x^1=x^2=x^3=0$ in \mythetag{5.12} and \mythetag{5.23},
we obtain
$$
\xalignat 3
&\hskip -2em
\Vmatrix X_0\\X_1\\X_2\\ X_3\endVmatrix
=\Vmatrix 0\\0\\0\\0\endVmatrix,
&&\Vmatrix Y_{01}\\Y_{02}\\Y_{03}\endVmatrix
=\Vmatrix 0\\0\\0\endVmatrix,
&&\Vmatrix Y_{12}\\Y_{13}\\Y_{23}\endVmatrix
=\Vmatrix 0\\0\\-4\,R^2\endVmatrix.
\quad
\mytag{5.24}
\endxalignat
$$
The quantities \mythetag{5.24} are initial data for the equations
\mythetag{3.20} and \mythetag{3.22}.\par
\head
6. Analysis of the compatibility conditions.
\endhead
     Six linearly independent Killing vector fields \mythetag{5.1},
\mythetag{5.5}, \mythetag{5.9}, \mythetag{5.13}, \mythetag{5.17},
and \mythetag{5.21} do always exist regardless to the function
\mythetag{1.4}. However, all of them are space-like vector fields
since $X_0=0$ for them. It is known that time-like Killing vector 
fields are more important for quantum field theories. For this reason
we look for the solutions of the equations \mythetag{4.8} and 
\mythetag{4.9} with
$$
\hskip -2em
X_0\neq 0.
\mytag{6.1}
$$
Under the assumption \mythetag{6.1} the second equation in 
\mythetag{4.8} produces two mutually exclusive options for the
function $R=R(x^0)$ in \mythetag{1.4}:
$$
\hskip -2em
R'=0\qquad\text{or}\qquad 2\,(R')^2-R''\,R+R^2=0.
\mytag{6.2}
$$
We study these two options in \mythetag{6.2} as two different cases.
\par
    {\bf The first case}: $R'=0$. In this case $R=\const$ and $R>0$,
i\.\,e\. $R$ is a positive constant. Applying the condition $R'=0$
to \mythetag{4.8}, we find that both of the equations \mythetag{4.8} 
are fulfilled identically in this case. As for \mythetag{4.9},
here we have
$$
\hskip -2em
2\,(R')^2-R''\,R+R^2=R^2\neq 0.
\mytag{6.3}
$$
Due to \mythetag{6.3} and due to $R'=0$ from \mythetag{4.9} we
derive
$$
\xalignat 3
&\hskip -2em
Y_{01}=0, &&Y_{02}=0, &&Y_{03}=0.
\mytag{6.4}
\endxalignat
$$
Now let's write the differential equations \mythetag{3.22}, taking 
into account \mythetag{6.4}. For the function $X_0$ we get the
following equations:
$$
\xalignat 4
&\frac{\partial X_0}{\partial x^0}=0,
&&\frac{\partial X_0}{\partial x^1}=0,
&&\frac{\partial X_0}{\partial x^2}=0,
&&\frac{\partial X_0}{\partial x^3}=0.
\qquad
\mytag{6.5}
\endxalignat
$$
The equations \mythetag{6.5} mean that $X_0$ is a constant function:
$$
X_0=\const.
$$
Moreover, taking into account \mythetag{6.4}, from \mythetag{3.20} 
and \mythetag{3.22} we derive:
$$
\xalignat 3
&\frac{\partial X_1}{\partial x^0}=0,
&&\frac{\partial X_2}{\partial x^0}=0,
&&\frac{\partial X_3}{\partial x^0}=0,
\qquad
\mytag{6.6}\\
\vspace{2ex}
&\frac{\partial Y_{12}}{\partial x^0}=0,
&&\frac{\partial Y_{23}}{\partial x^0}=0,
&&\frac{\partial Y_{13}}{\partial x^0}=0.
\qquad
\mytag{6.7}
\endxalignat
$$
The equations \mythetag{6.6} and \mythetag{6.7} mean that the 
functions $X_1,\,X_2,\,X_3,\,Y_{12},\,Y_{13},\,Y_{23}$ actually 
do not depend on the variable $x^0$. As appears, other equations
in \mythetag{3.20} and \mythetag{3.22} in the case of $R'=0$ do 
not contain $X_0$ and form complete system of Pfaff equations 
for six functions $X_1,\,X_2,\,X_3,\,Y_{12},\,Y_{13},\,Y_{23}$
with respect to three variables $x^1,\,x^2,\,x^3$. These equations
have at most six linearly independent solutions. These solutions 
are exhausted by six Killing vector field considered in 
section~5.
\mytheorem{6.1} In the case of $R'=0$ the spherical universe
$M=\Bbb R\times S^3$ with the metric \mythetag{1.2} admits 
exactly one linearly independent time-like Killing vector 
field 
$$
\hskip -2em
\bold X=\frac{\partial}{\partial x^0}
\mytag{6.8}
$$
in addition to six space-like Killing vector fields \mythetag{5.1},
\mythetag{5.5}, \mythetag{5.9}, \mythetag{5.13}, \mythetag{5.17}, 
\mythetag{5.21} produced by the rotations of the sphere $S^3$. 
\endproclaim
The vector field \mythetag{6.8} is orthogonal to the sphere $S^3$
in $M$. It commutes with other six Killing vector fields, which are 
tangent to $S^3$.\par
    {\bf The second case}. In contrast to \mythetag{6.3}, in this
case we have the following equality for the function $R=R(x^0)$ in 
\mythetag{1.4}:
$$
\pagebreak
\hskip -2em
2\,(R')^2-R''\,R+R^2=0.
\mytag{6.9}
$$
Due to the equality \mythetag{6.9} the compatibility conditions
\mythetag{4.9} and the first compatibility equation \mythetag{4.8}
are fulfilled identically. Moreover, we have 
$$
\hskip -2em
\gathered
-(2\,(R')^2-R''\,R+R^2)'+2\,R'\,(2\,(R')^2-R''\,R+R^2)=\\
=4\,(R')^3-5\,R''\,R'\,R+\,R'''\,R^2.
\endgathered
\mytag{6.10}
$$
Due to \mythetag{6.10} the first compatibility equation in 
\mythetag{4.8} is also fulfilled identically. Thus, the 
equation \mythetag{6.9} is the only compatibility condition
derived from \mythetag{3.27} in the second case.\par
    Note that the equation \mythetag{6.9} can be integrated 
up to the first order differential equation. Indeed, since
$R\neq 0$, it can be written as follows:
$$
\hskip -2em
\left(\frac{R'}{R^2}\right)\kern -4pt\raise 6.5pt\hbox{$'$}
=\frac{1}{R}
\mytag{6.11}
$$
Let's multiply both sides of \mythetag{6.11} by the fraction
$$
\hskip -2em
\frac{R'}{R^2}
\mytag{6.12}
$$
As a result we get the equation with the pure derivatives in both 
sides:
$$
\hskip -2em
\left(\frac{1}{2}\left(\frac{R'}{R^2}\right)\kern -2pt\raise 9.5pt
\hbox{$\ssize 2$}\right)\kern -4pt\raise 6.5pt\hbox{$'$}=
\frac{R'}{R^3}=\left(-\frac{1}{2}\,\frac{1}{R^2}\right)\kern -4pt
\raise 6.5pt\hbox{$'$}.
\mytag{6.13}
$$
Integrating the equality \mythetag{6.13}, we derive
$$
\hskip -2em
(R')^2=C\,R^4-R^2,
\mytag{6.14}
$$
where $C$ is a constant of integration. Note that $R'$ can vanish 
at some points, it is not identically zero in this case. The same 
is true for the fraction \mythetag{6.12}. Therefore the equation 
\mythetag{6.14} is equivalent to the initial equation \mythetag{6.9} 
at all point except for those, where $R'=0$.\par
    It is clear that $C$ in \mythetag{6.14} is a positive constant.
Let's denote $C=1/a^2$. Then $R$ is a function with the values ranging
in the interval
$$
R\in [a,+\infty).
$$
The equation \mythetag{6.14} itself can be written as follows
$$
\hskip -2em
\left(\left(\frac{1}{R}\right)\kern -4pt\raise 6.5pt\hbox{$'$}
\right)\kern -2pt\raise 9.5pt\hbox{$\ssize 2$}=\frac{1}{a^2}
-\frac{1}{R^2}.
\mytag{6.15}
$$
Let's denote $u=1/R$ for a while. Then we transform \mythetag{6.15}
to
$$
\hskip -2em
(u')^2=1/a^2-u^2
\mytag{6.16}
$$
The equation \mythetag{6.16} can be integrated. Its general solution
looks like
$$
\hskip -2em
u(x^0)=\frac{\cos(x^0+b)}{a},
\mytag{6.17}
$$
where $b$ is a constant of integration. Without loss of generality
we can take $b=0$. Then from \mythetag{6.17} we derive the following
formula
$$
\hskip -2em
R(x^0)=\frac{a}{\cos(x^0)}.
\mytag{6.18}
$$\par
    Having defined the function \mythetag{1.4} by means of the formula
\mythetag{6.18}, now let's define the time variable $t$ by means of
the differential equation \mythetag{1.5}:
$$
\hskip -2em
\frac{dx^0}{\cos(x^0)}=\frac{c\,dt}{a}.
\mytag{6.19}
$$
Integrating both sides of the equality \mythetag{6.19}, we obtain:
$$
\hskip -2em
\ln\left(\frac{1+\sin(x^0)}{\cos(x^0)}\right)=\frac{c\,t}{a}.
\mytag{6.20}
$$
Transforming \mythetag{6.20}, we pass from logarithms to exponentials.
As a result we get:
$$
\hskip -2em
\frac{1+\sin(x^0)}{\cos(x^0)}=e^{\raise 2pt\hbox{$\ssize\frac{c\,t}{a}$}}.
\mytag{6.21}
$$
Now we square both sides of the equality \mythetag{6.21}. This yields
$$
\hskip -2em
\frac{1+2\,\sin(x^0)+\sin^2(x^0)}{\cos^2(x^0)}
=\frac{2+2\,\sin(x^0)-\cos^2(x^0)}{\cos^2(x^0)}
=e^{\raise 2pt\hbox{$\ssize\frac{2\,c\,t}{a}$}}.
\mytag{6.22}
$$
The equality \mythetag{6.22} can be transformed to the following one:
$$
\hskip -2em
1+\sin(x^0)=\frac{1+e^{\raise 2pt
\hbox{$\ssize\frac{2\,c\,t}{a}$}}}{2}\,\cos^2(x^0).
\mytag{6.23}
$$
Note that the left hand side of \mythetag{6.23} coincides with
the numerator of the fraction in the left hand side of \mythetag{6.21}.
Substituting \mythetag{6.23} back into \mythetag{6.21}, we get
$$
\hskip -2em
\cos(x^0)=\frac{2\,e^{\raise 2pt\hbox{$\ssize\frac{c\,t}{a}$}}}
{1+e^{\raise 2pt\hbox{$\ssize\frac{2\,c\,t}{a}$}}\vphantom{\vrule
height 13pt}}=\frac{1}{\cosh\bigl(\frac{c\,t}{a}\bigr)
\vphantom{\vrule height 13pt}}.
\mytag{6.24}
$$
Substituting \mythetag{6.24} into \mythetag{6.18} we find the
dependence of $R$ on the time variable $t$:
$$
R(t)=a\,\cosh\bigl(\tsize\frac{c\,t}{a}\bigr).
$$
Moreover, substituting \mythetag{6.24} into \mythetag{6.21},
we derive the following formula:
$$
\hskip -2em
\sin(x^0)=\frac{\sinh\bigl(\frac{c\,t}{a}\bigr)\vphantom{\vrule 
depth 9pt height 0pt}}{\cosh\bigl(\frac{c\,t}{a}\bigr)\vphantom{\vrule 
height 13pt}}=\tanh\bigl(\tfrac{c\,t}{a}\bigr).
\mytag{6.25}
$$\par 
    Now, relying on the above calculations, we introduce the modified
stereographic coordinates especially for this particular case:
$$
\xalignat 4
&\hskip -2em
u^0=\frac{c\,t}{a}, 
&&u^1=x^1,
&&u^2=x^2,
&&u^3=x^3.
\qquad
\mytag{6.26}
\endxalignat
$$
The metric \mythetag{1.2} in these coordinates $u^0,\,u^1,\,u^2,\,u^3$
is written as follows:
$$
ds^{\kern 0.5pt 2}=a^2\,(du^0)^2
-4\,a^2\cosh^2(u^0)\,\frac{(du^1)^2+(du^2)^2+(du^3)^2}
{\left(|u|^2+1\right)^2},
\mytag{6.27}
$$
where $|u|^2=(u^1)^2+(u^2)^2+(u^3)^2$. Let's consider the 
following four functions of the modified stereographic coordinates 
$u^0,\,u^1,\,u^2,\,u^3$:
$$
\xalignat 2
&\hskip -2em
z^1=A\,\cosh(u^0)\,\frac{2\,u^1}{|u|+1},
&&z^2=A\,\cosh(u^0)\,\frac{2\,u^2}{|u|+1},
\qquad\\
\vspace{-1ex}
&&&\mytag{6.28}\\
\vspace{-1ex}
&\hskip -2em
z^3=A\,\cosh(u^0)\,\frac{2\,u^3}{|u|+1},
&&z^4=A\,\cosh(u^0)\,\frac{|u|-1}{|u|+1}.
\qquad
\endxalignat
$$
As it was shown in \mycite{4}, the functions \mythetag{6.28} determine
an embedding of the sphere $S^3$ into the four-dimensional Euclidean
space $\Bbb R^4$. Let's complement the functions \mythetag{6.28} with 
one additional function
$$
\hskip -2em
z^0=A\,\sinh(u^0).
\mytag{6.29}
$$
Then the functions \mythetag{6.28} and \mythetag{6.29} taken together
determine an embedding of our universe $M=\Bbb R\times S^3$ into the
five-dimensional space $\Bbb R^5$. If we equip this space with the
sign-indefinite metric
$$
\hskip -2em
ds^2=(dz^0)^2-(dz^1)^2-(dz^2)^2-(dz^3)^2-(dz^4)^2,
\mytag{6.30}
$$
then we find that the metric \mythetag{6.30} induces the metric
\mythetag{6.27} in $M$ under the embedding given by the functions
\mythetag{6.28} and \mythetag{6.29}.\par
     Calculating the curvature tensor for the metric \mythetag{6.27},
we find that our universe $M=\Bbb R\times S^3$ in this case is a
manifold of constant negative sectional curvature 
$$
K=-\frac{1}{a^2}.
$$
It is known that any four-dimensional constant curvature manifold
has exactly ten linear independent Killing vector fields (compare
this fact with the theorem~\mythetheorem{3.1}). Six of them are given 
by the formulas \mythetag{5.1}, \mythetag{5.5}, \mythetag{5.9},
\mythetag{5.13}, \mythetag{5.17}, \mythetag{5.21}. These fields are
associated with meridional and equatorial rotations of the sphere 
$S^3$. Other four fields are determined by the hyperbolic rotations
of $M$ itself.\par
    {\bf Hyperbolic rotation} in the plane $z^1 O z^0$. This rotation 
induces the Killing vector field $\bold X$ in $M$ expressed by the 
formula
$$
\hskip -2em
\gathered
\bold X=\frac{2\,(u^1)}{|u|^2+1}\,\frac{\partial}{\partial u^0}+
\frac{|u|^2-2\,(u^1)^2+1}{2}\,\frac{\sinh(u^0)}{\cosh(u^0)}
\,\frac{\partial}{\partial u^1}\,-\\
\vspace{1ex}
-\,(u^1)(u^2)\,\frac{\sinh(u^0)}{\cosh(u^0)}
\,\frac{\partial}{\partial u^2}
-\,(u^1)(u^3)\,\frac{\sinh(u^0)}{\cosh(u^0)}
\,\frac{\partial}{\partial u^3}.
\endgathered
\mytag{6.31}
$$\par
    {\bf Hyperbolic rotation} in the plane $z^2 O z^0$. This rotation 
induces a Killing vector field in $M$ very similar to the previous one. 
\pagebreak The Killing vector field $\bold X$ for this case is expressed 
by the following formula:
$$
\hskip -2em
\gathered
\bold X=\frac{2\,(u^2)}{|u|^2+1}\,\frac{\partial}{\partial u^0}
-\,(u^2)(u^1)\,\frac{\sinh(u^0)}{\cosh(u^0)}
\,\frac{\partial}{\partial u^1}\,+\\
\vspace{1ex}
+\,\frac{|u|^2-2\,(u^2)^2+1}{2}\,\frac{\sinh(u^0)}{\cosh(u^0)}
\,\frac{\partial}{\partial u^2}
-(u^2)(u^3)\,\frac{\sinh(u^0)}{\cosh(u^0)}
\,\frac{\partial}{\partial u^3}.
\endgathered
\mytag{6.32}
$$\par
    {\bf Hyperbolic rotation} in the plane $z^3 O z^0$. This rotation 
induces the Killing vector field $\bold X$ in $M$ expressed by the 
formula
$$
\hskip -2em
\gathered
\bold X=\frac{2\,(u^3)}{|u|^2+1}\,\frac{\partial}{\partial u^0}
-\,(u^3)(u^1)\,\frac{\sinh(u^0)}{\cosh(u^0)}
\,\frac{\partial}{\partial u^1}\,-\\
\vspace{1ex}
-\,(u^3)(u^2)\,\frac{\sinh(u^0)}{\cosh(u^0)}
\,\frac{\partial}{\partial u^2}
+\frac{|u|^2-2\,(u^3)^2+1}{2}\,\frac{\sinh(u^0)}{\cosh(u^0)}
\,\frac{\partial}{\partial u^3}.
\endgathered
\mytag{6.33}
$$\par
    {\bf Hyperbolic rotation} in the plane $z^4 O z^0$. This rotation 
produces the Killing vector field $\bold X$ in $M$ expressed by the 
formula
$$
\hskip -2em
\gathered
\bold X=\frac{|u|^2-1}{|u|^2+1}\,\frac{\partial}{\partial u^0}
+(u^1)\,\frac{\sinh(u^0)}{\cosh(u^0)}
\,\frac{\partial}{\partial u^1}\,+\\
\vspace{1ex}
+\,(u^2)\,\frac{\sinh(u^0)}{\cosh(u^0)}
\,\frac{\partial}{\partial u^2}
+(u^3)\,\frac{\sinh(u^0)}{\cosh(u^0)}
\,\frac{\partial}{\partial u^3}.
\endgathered
\mytag{6.34}
$$\par
Using \mythetag{6.26}, \mythetag{6.20}, \mythetag{6.24}, and 
\mythetag{6.25}, one can easily transform the above four vector
fields to the initial North Pole stereographic coordinates
$x^0,\,x^1,\,x^2,\,x^3$. Note that none of the vector fields
\mythetag{6.31}, \mythetag{6.32}, \mythetag{6.33}, \mythetag{6.34}
is purely time-like. They are build by time-like vectors at some
points of $M$ and by space-like vectors at some other points.
\mytheorem{6.2} In the case of $2\,(R')^2-R''\,R+R^2=0$ the spherical
universe $M=\Bbb R\times S^3$ with the metric \mythetag{1.2} admits 
four linearly independent Killing vector fields 
\mythetag{6.31}, \mythetag{6.32}, \mythetag{6.33}, \mythetag{6.34}
in addition to six space-like Killing vector fields \mythetag{5.1},
\mythetag{5.5}, \mythetag{5.9}, \mythetag{5.13}, \mythetag{5.17}, 
\mythetag{5.21}. Neither of these ten Killing vector fields, nor any 
linear combination of them with constant coefficients is a purely 
time-like vector field in $M$.
\endproclaim
\head
7. Conclusions.
\endhead
     The main result of this paper is that in general case a 
homogeneous and isotropic closed universe $M=\Bbb R\times S^3$ 
has no time-like Killing vector fields at all, i\.\,e\. it is 
non-stationary in the sense of \mycite{1} and \mycite{2}. For 
this reason it is a good model for to study various quantization 
procedures for the matter fields in the presence of a non-stationary 
gravitation field as a background. The theorems~\mythetheorem{6.1} 
and \mythetheorem{6.2} specify two exceptional cases. In the first of 
them the universe $M=\Bbb R\times S^3$ is stationary in whole, 
while in the second case is is piecewise stationary, i\.\,e\.
it is broken into stationary and non-stationary fragments.

\enddocument
\end